\documentclass[12pt]{article}
\usepackage{amsmath,amssymb,amsfonts,latexsym,amssymb}
\usepackage[dvips]{epsfig}
\usepackage{amscd}
\usepackage{enumerate}

\title{A brief review of abelian categorifications}
\author{Mikhail Khovanov,
Volodymyr Mazorchuk\\ and Catharina Stroppel}

\date{}

\newtheorem{prop}{Proposition}

\newtheorem{theorem}[prop]{Theorem}
\newtheorem{definition}[prop]{Definition}

\newcommand{\oplusop}[1]{{\mathop{\oplus}\limits_{#1}}}
\newcommand{\oplusoop}[2]{{\mathop{\oplus}\limits_{#1}^{#2}}}

\begin{document}

\maketitle

\def\C{\mathbb C}
\def\R{\mathbb R}
\def\N{\mathbb N}
\def\Z{\mathbb Z}
\def\Q{\mathbb Q}
\def\g{\mathfrak g}
\def\p{\mathfrak p}
\def\h{\mathfrak h}
\def\n{\mathfrak n}
\newcommand{\Ext}{\operatorname{Ext}}
\newcommand{\End}{\operatorname{End}}
\newcommand{\add}{\operatorname{add}}

\def\F{\mathbb F}
\def\S{\mathbb S}
\def\l{\lbrace}
\def\r{\rbrace}
\def\o{\otimes}
\def\lra{\longrightarrow}
\newcommand{\ba}{\mathbf{a}}
\newcommand{\dmod}{\mathrm{-mod}}
\newcommand{\mc}{\mathcal}
\def\cF{\mathcal{F}}
\def\Hom{\textrm{Hom}}
\def\drawing#1{\begin{center} \epsfig{file=#1} \end{center}}
\def\mc{\mathcal}
\def\mf{\mathfrak}
\def\mb{\mathbb}

\def\yesnocases#1#2#3#4{\left\{ \begin{array}{ll} #1 & #2 \\ #3 & #4
\end{array} \right. }

\newcommand{\define}{\stackrel{\mbox{\scriptsize{def}}}{=}}
\def\hsm{\hspace{0.05in}}

\def\cO{\mc{O}}   
\def\cC{\mc{C}}
\def\sln{\mathfrak{sl}(n)}

\begin{abstract}
This article contains a review of categorifications of semisimple
representations of various rings via abelian categories and exact
endofunctors on them. A simple definition of an abelian
categorification is presented and illustrated with several examples,
including categorifications of various representations of the
symmetric group and its Hecke algebra via highest weight categories
of modules over the Lie algebra $\mathfrak{sl}_n$. The review is
intended to give non-experts in representation theory who are
familiar with the topological aspects of categorification (lifting
quantum link invariants to homology theories) an idea for the sort
of categories that appear when link homology is extended to tangles.
\end{abstract}




\section{A simple framework for categorification} \label{sec-framework}

{\bf Categorification.} The Grothendieck group $K(\mc{B})$ of an
abelian category $\mc{B}$ has as generators the symbols $[M],$ where
$M$ runs over all the objects of $\mc{B},$ and defining relations
$[M_2]=[M_1]+ [M_3],$ whenever there is a short exact sequence
$$ 0 \lra M_ 1\lra M_2 \lra M_3 \lra 0.$$
An exact functor $F$ between abelian categories induces a homomorphism $[F]$
between their Grothendieck groups.

Let $A$ be a ring which is free as an abelian group, and $\ba=\{
a_i\}_{i\in I}$ a basis of $A$ such that the multiplication has
nonnegative integer coefficients in this basis:
\begin{eqnarray}
\label{cs} a_i a_j = \sum_k c_{ij}^k a_k, \hspace{0.2in} c_{ij}^k
\in \Z_{\geq0}.
\end{eqnarray}
Let $B$ be a (left) $A$-module.

\begin{definition} \label{def-1}
A (weak) abelian categorification of $(A,\ba,B)$ consists of an
abelian category $\mc{B},$  an isomorphism $\varphi:K(\mc{B})
\overset{\sim}{\longrightarrow}B$ and exact endofunctors $F_i:
\mc{B} \lra \mc{B}$ such that the following holds:
\begin{enumerate}[(C-I)]
\item \label{I}
The functor $F_i$ lifts the action of $a_i$ on the module $B$, i.e. the action
of $[F_i]$ on the Grothendieck group of $\mc{B}$ descends to the
action of $a_i$ on the module $B$ so that the diagram below is
commutative.
 $$\begin{CD}
    K(\mc{B})  @>{[F_i]}>>     K(\mc{B})   \\
   @V{\varphi}VV                   @VV{\varphi}V   \\
    B   @>{a_i}>>     B
  \end{CD} $$
\item \label{II}
There are isomorphisms
$$ F_i F_j \cong \oplusop{k} F_k^{c_{ij}^k},$$
i.e., the composition $F_i F_j$  decomposes as the direct sum of
functors $F_k$ with multiplicities $c_{ij}^k$ as in \eqref{cs}.
\end{enumerate}
\end{definition}

If there is a categorification as above we say {\it the action of
the functors $F_i$ on the category $\mc{B}$ categorifies the
$A$-module $B$}.\\

In all our examples, the objects of $\mc{B}$ will have finite length
(finite Jordan-H\"older series). Consequently, if $\{L_j\}_{j\in J}$
is a collection of  simple objects of $\mc{B},$ one for each
isomorphism class, the  Grothendieck group $K(\mc{B})$ is free
abelian with basis elements $[L_j].$ The image of any object $M\in
\mc{B}$ in the Grothendieck group is
$$ [M] = \sum_j m_j(M) [L_j]$$
where $m_j(M)$ is the multiplicity of $L_j$ in some (and hence
in any) composition series of $M.$

The free group $K(\mc{B})$ has therefore a distinguished basis
$[L_j]_{j\in J}$, and the action of $[F_i]$ in this basis has
integer non-negative coefficients: $$[F_i (L_j)] =\sum
d_{ij}^k[L_k],$$ with $d_{ij}^k$ being the multiplicity
$m_k(F_i(L_j)).$ Via the isomorphism $\varphi$ we obtain a
distinguished basis $\mathbf{b}=\{b_j\}_{j\in J}$ of $B$, and
\begin{eqnarray}
\label{ab}
a_ib_j =\sum d_{ij}^k b_k.
\end{eqnarray}

Conversely, we could fix a basis $\mathbf{b}$ of $B$ with a
positivity constraint for the action of $A.$ as in \eqref{ab}. Then
our definition of a categorification of $(A,\ba,B)$ can be amended
to a similar definition of a categorification of
$(A,\ba,B,\mathbf{b}),$ with the additional data being the fixed
basis $\mathbf{b}$. Ideally the basis $\mathbf{b}$ corresponds then
via the isomorphism $\varphi$ to a basis $[M_j]_{j\in J}$ for
certain objects $M_j\in \mc{B}$. Varying the choice of basis might
give rise to an interesting combinatorial interplay between several,
maybe less prominent than $[L_j]_{j\in J}$ but more interesting,
families $\{M_j\}_{j\in J}$ of objects in $\mc{B}$. Typical examples
of such an interplay can be found in \cite{Brundan-bases}, \cite[Section 5]{FKS}.

\vspace{0.1in}

Of course, any such data $(A,\ba, B,\mathbf{b})$
admits a rather trivial categorification, via a semisimple category $\mc{B}.$
Namely, choose a field $\Bbbk$
and denote by $\Bbbk\mathrm{-vect}$ the category of finite-dimensional
$\Bbbk$-vector spaces. Let
$$\mc{B}= \bigoplus_{j\in J} \Bbbk\mathrm{-vect}$$
be the direct sum of categories $\Bbbk\mathrm{-vect},$ one for each
basis vector of $B.$ The category $\mc{B}$ is semisimple, with
simple objects $L_j$ enumerated by elements of $J,$ and
$$ \Hom_{\mc{B}}(L_j, L_k) =
{\left\{ \begin{array}{ll} \Bbbk &
\mathrm{if \hspace{0.1in}}j=k, \\ 0 & \mathrm{otherwise}.
\end{array}\right. } $$
We identify $K(\mc{B})$ with $B$ by mapping $[L_j]$ to $b_j$ The
functors $F_i$ are determined by their action on
simple objects, hence, given \eqref{ab}, we can define
$$ F_i(L_j) = \oplusop{k\in J} L_k^{d_{ij}^k}$$
and obtain a categorification of $(A,\ba, B,\mathbf{b})$. With few
exceptions, semisimple categorifications bring little or no new
structure into play, and we will ignore them. More interesting
instances of categorifications occur for non-semisimple categories
$\mc{B}.$ Here is a sample list.

\vspace{0.1in}


{\bf 1.} Let $\mathcal{A}_1$ be the first Weyl algebra (the algebra of polynomial differential
operators in one variable) with integer coefficients,
$$\mathcal{A}_1 = \Z  \langle x, \partial
\rangle /(\partial x - x\partial -1).$$ We fix the basis $\{x^i
\partial^j\}_{i,j\ge 0}$ of $\mathcal{A}_1$. The $\Z$-lattice
$B\subset \Q[x]$ with the basis $\mathbf{b}=\{\frac{x^n}{n!}\}_{n\ge
0}$ is an $\mathcal{A}_1$-module.

To categorify this data we consider the category $\mc{B}= \oplusop{n
\ge 0} R_n\mathrm{-mod}$, i.e. the direct sum of the categories of
finite-dimensional modules over the nilCoxeter $\Bbbk$-algebra
$R_n.$ The latter has generators $Y_1, \dots, Y_{n-1}$ subject to
relations
\begin{eqnarray*}
    Y_i^2 & = & 0,   \\
    Y_i Y_j & = & Y_j Y_i \hspace{0.2in} \mathrm{if} \hspace{0.2in}
   |i-j|>1, \\
    Y_i Y_{i+1} Y_i & = & Y_{i+1} Y_i Y_{i+1}.
\end{eqnarray*}
The algebra $R_n$ has a unique, up to isomorphism, finite
dimensional simple module $L_n,$ and $K(R_n\mathrm{-mod})\cong \Z.$
The Grothendieck group $K(\mc{B})$ is naturally isomorphic to the
$\mathcal{A}_1$-module $B,$ via the isomorphism $\varphi$ which maps
$[L_n]$ to $\frac{x^n}{n!}.$ The endofunctors $X,D$ in $\mc{B}$ that
lift the action of $x$ and $\partial$ on $B$ are the induction and
restriction functors for the inclusion of algebras $R_n \subset
R_{n+1}.$ One takes $\{x^i \partial^j\}_{i,j\ge 0} $ as the
basis $\ba$ of $\mathcal{A}_1.$ Basis elements lift to functors $X^i D^j$,
and the isomorphisms (C-II) of definition~\ref{def-1} are induced by
an isomorphism of functors $ D X \cong X D\oplus
\mathrm{Id}$ which lifts the defining relation $\partial x = x \partial + 1$ in
 $\mathcal{A}_1.$ A detailed analysis of this categorification can be found in
\cite{Kh-nilcoxeter}.

 \vspace{0.15in}


{\bf 2.} The regular representation of the group ring $\Z[S_n]$ of
the symmetric group $S_n$ has a categorification via projective
functors acting on a regular block of the highest weight BGG
category $\mc{O}$ from \cite{BGG} for $\mf{sl}_n$. (For an
introduction to the representation theory of semisimple Lie algebras
we refer
to \cite{Humphreys}).\\

To define category $\mc{O}$, start with the standard triangular
decomposition $\mf{sl}_n= \mf{n}_+\oplus \mf{h}\oplus \mf{n}_-,$
where the first and the last terms are the Lie algebras of strictly
upper-triangular (resp. lower-triangular) matrices, while $\mf{h}$
is the algebra of traceless diagonal matrices. The highest weight
category $\mc{O}$ of $\mf{sl}_n$ is the full subcategory of the
category of finitely-generated $\mf{sl}_n$-modules consisting of
$\mf{h}$-diagonalisable (possibly infinite dimensional) modules on
which $U(\mf{n}_+)$ acts locally-nilpotently. Thus, any $M\in
\mc{O}$ decomposes as
$$M=\oplusop{\lambda\in \mf{h}^{\ast}} M_{\lambda},$$
where $hx= \lambda(h)x$ for any $h\in \mf{h}$ and $x\in M_{\lambda}.$
Here $\mf{h}^{\ast}$ is the dual vector space of $\mf{h}$, its elements
are called weights.

The one-dimensional modules $\underline{\C}_{\lambda}=\C
v_{\lambda}$ over the positive Borel subalgebra $\mf{b}_+=
\mf{n}_+\oplus \mf{h}$ are classified by elements $\lambda$ of
$\mf{h}^{\ast}.$ The subalgebra $\mf{n}_+$ acts trivially on
$v_{\lambda},$ while $hv_{\lambda} = \lambda(h)v_{\lambda}$ for
$h\in \mf{h}.$

The {\it Verma module} $M(\lambda)$ is the $\mf{sl}_n$-module
induced from the $\mf{b}_+$-module $\underline{\C}_{\lambda},$
$$M(\lambda) = U(\mf{sl}_n) \otimes_{U(\mf{b}_+)}
 \underline{\C}_{\lambda}.$$
The Verma module $M(\lambda)$  has a unique simple quotient, denoted
$L(\lambda),$ and any simple object of $\mc{O}$ is isomorphic to
$L(\lambda)$ for some $\lambda\in \mf{h}^{\ast}.$

We call a weight $\lambda$ \emph{positive integral} if $\langle
\lambda, \alpha\rangle \in \Z_{\geq0}$ for any positive simple root
$\alpha\in\mathfrak{h}^*.$ The representation $L(\lambda)$ is
finite-dimensional if and only if $\lambda$ is a positive integral
weight.

Although most of the objects in $\cO$ are infinite dimensional
vector spaces,  every object $M$ of $\mc{O}$ has finite length, i.e.
there is an increasing filtration by subobjects $0=M^0 \subset
M^1\subset \dots \subset M^m = M$ such that the subsequent quotients
$M^{i+1}/M^i$ are isomorphic to simple objects, hence have the form
$L(\lambda)$ (where $\lambda$ may vary). The Grothendieck group of
$\mc{O}$ is thus a free abelian group with generators $[L(\lambda)]$
for $\lambda\in \mf{h}^{\ast}.$

It turns out that $\mc{O}$ has enough projective objects: given $M$
there exists a surjection $P \twoheadrightarrow M$ with a projective
$P\in\mc{O}$. Moreover, isomorphism classes of indecomposable
projective objects are enumerated by elements of $\mf{h}^{\ast}.$
The indecomposable projective object $P(\lambda)$ is determined by
the property of being projective and
\begin{equation*}
\mathrm{Hom}_{\mc{O}}(P(\lambda), L(\mu)) =
\begin{cases} \C & \text{if $\lambda=\mu$,} \\
0 &\text{otherwise.}
\end{cases}
\end{equation*}
We should warn the reader that the $P(\lambda)$'s are not projective
when viewed as objects of the category of {\it all}
$\mf{sl}_n$-modules, while the $L(\lambda)$'s remain simple in the
latter category.

The symmetric group $S_n,$ the Weyl group of $\mf{sl}_n,$ acts naturally
on $\mf{h}$ by permuting the diagonal entries and then also on
$\mf{h}^*$. Let $\rho\in \mf{h}^{\ast}$ be the half-sum of positive roots.
In the study of the category $\mc{O}$ an important role is played by
the shifted (dot) action of $S_n,$
$$ w\cdot \lambda = w(\lambda+\rho) -\rho.$$
Two simple modules $L(\lambda), L(\lambda')$ have the same central
character (i.e. are annihilated by the same maximal ideal of the
centre of the universal enveloping algebra) if and only if $\lambda$
and $\lambda'$ belong to the same $S_n$-orbit under the shifted
action. Consequently, $\mc{O}$ decomposes into a direct sum of
categories
\begin{equation} \label{dir-eq}
\mc{O} = \oplusop{\nu \in \mf{h}^{\ast}/S_n} \mc{O}_{\nu}
\end{equation}
indexed by orbits $\nu$ of the shifted action of $S_n$ on
$\mf{h}^{\ast}$. Here, $\mc{O}_{\nu}$ consists of all modules with
composition series having only simple subquotients isomorphic to
$L(\lambda)$ for $\lambda \in \nu.$ There is no interaction between
$\mc{O}_{\nu}$ and  $\mc{O}_{\nu'}$ for different orbits $\nu,
\nu'.$ More accurately, if $\nu\not=\nu'$ then $\mathrm{Ext}_{\mc
O}^i(M,M')=0$ for any $i\geq 0$, $M\in \mc{O}_{\nu}$ and $M'\in
\mc{O}_{\nu'}$.

Furthermore, each $\mc{O}_{\nu}$ is equivalent to the category of
finite-dimensional modules over some finite-dimensional $\C$-algebra
$A_{\nu}.$ Here's the catch, though: {\it explicitly} describing
$A_{\nu}$ for $n>3$ and interesting $\nu$ is very hard, see
\cite{St-quiver}. For an {\it implicit} description we just form
$P=\oplusop{\lambda\in \nu}P(\lambda),$ the direct sum of all
indecomposable projectives over $\lambda\in \nu.$ Then $A_{\nu}\cong
\mathrm{Hom}_{\mc O}(P,P)^{op}.$

An orbit $\nu$ (for the shifted action) is called \emph{generic} if
$w\cdot \lambda-\lambda$ is never integral, for $\lambda \in \nu$
and $w\in S_n, w\not= 1.$ For a generic orbit $\nu,$ the category
$\mc{O}_{\nu}$ is boring and equivalent to the direct sum of $n!$
copies of the category of finite-dimensional $\C$-vector spaces, one
for each $\lambda\in \nu.$ For such $\lambda$ we have $P(\lambda) =
M(\lambda) = L(\lambda),$ i.e. the Verma module with the highest
weight $\lambda$ is simple as well as projective in $\mc{O}.$

We call an orbit \emph{integral} if it is a subset of the integral
weight lattice in $\mf{h}^{\ast}.$  In \cite{Soergel-garben} it is
shown that $\mc{O}_{\nu}$ for non-integral $\nu$ reduces
to those for integral weights. From now on we
therefore assume that $\nu$ is integral. Then the category
$\mc{O}_{\nu}$ is indecomposable (unlike in the generic case).
Moreover, the complexity of $\mc{O}_{\nu}$ only depends on the type
of the orbit. If two orbits $\nu$ and $\nu'$ contain points
$\lambda\in\nu$, $\lambda'\in\nu'$ with identical stabilisers, then
the categories $\mc{O}_{\nu}$ and $\mc{O}_{\nu'}$ are equivalent,
see \cite{BeGe}, \cite{Soergel-garben}. If the stabiliser of $\nu$ under the
shifted action is trivial, the category $\mc{O}_{\nu}$ is called a
{\it regular block}. Regular blocks are the most complicated
indecomposable direct summands of $\mc{O},$ for instance in the
sense of having the maximal number of isomorphism classes of simple
modules.

There is a natural bijection between the following three sets:
positive integral weights, isomorphism classes of irreducible
finite-dimensional representations of $\mf{sl}_n,$ and regular
blocks of $\mc{O}$ for $\mf{sl}_n.$ A positive integral weight
$\lambda$ is the highest weight of an irreducible finite-dimensional
representation $L(\lambda),$ determined by the weight uniquely up to
isomorphism. In turn, $L(\lambda)$ belongs to the regular block
$\mc{O}_{\nu},$ where $\nu=S_n \cdot \lambda$ is the orbit of
$\lambda.$

Any two regular blocks of $\mc{O}$ are equivalent as categories,
as shown in~\cite{Jan}. For this reason, we can restrict our
discussion to the uniquely defined regular block which contains the
one-dimensional trivial representation $L(0)$ of $\mf{sl}_n.$ We
denote this block by $\mc{O}_0.$ It has $n!$ simple modules
$L(w)=L(w\cdot 0),$ enumerated by all permutations $w\in S_n$ (with
the identity element $e$ of $S_n$ corresponding to $L(0)$ which is
the only finite dimensional simple module in $\mc{O}_0$). Thus,
$K(\mc{O}_0)$ is free abelian of rank $n!$ with basis $\{
[L(w)]\}_{w\in S_n}.$ Other notable objects in $\mc{O}_0$ are the
Verma modules $M(w)=M(w\cdot 0)$ and the indecomposable projective
modules $P(w)=P(w\cdot 0),$ over all $w\in S_n.$ The sets
$\{[M(w)]\}_{w\in S_n}$ and $\{[P(w)]\}_{w\in S_n}$ form two other
prominent bases in $K(\mc{O}_0).$ For the set $\{[M(w)]\}_{w\in
S_n}$ this is easy to see, because the transformation matrix between
Verma modules and simple modules is upper triangular with ones on
the diagonal. For the set $\{[P(w)]\}_{w\in S_n}$ this claim is not
obvious and  relates to the fact that $\mc{O}_0$ has finite
homological dimension, see \cite{BGG}.

Equivalences between regular blocks are established by means of
translation functors. First note that we can tensor two
$U(\mf{sl}_n)$-modules over the ground field. If $V$ is a {\em
finite-dimensional} $\mf{sl}_n$-module it follows from the
definitions that $V\otimes M$ lies in $\mc{O}$ whenever $M$ is in
$\mc{O}.$ Hence, tensoring with $V$ defines an endofunctor
$V\otimes{}_-$ of the category $\mc{O}$. Taking direct summands of
the functors $V\otimes{}_-$ provides a bewildering collection of
different functors and allows one to analyse $\mc{O}$ quite deeply.
By definition, a {\it projective functor} is any endofunctor of
$\mc{O}$ isomorphic  to a direct summand of $V\otimes{}_-$ for some
finite-dimensional $\mf{sl}_n$-module $V.$ Projective functors were
classified by J.~Bernstein and S.~Gelfand \cite{BeGe}. {\it
Translation functors} are special cases of projective functors.

Let us restrict our discussion to projective endofunctors in the
regular block $\mc{O}_0.$  Each projective endofunctor $\mc{O}_0\lra
\mc{O}_0$ decomposes into a finite direct  sum of indecomposable
functors $\theta_w,$ enumerated by permutations $w$ and determined
by the property $\theta_w(M(e))\cong P(w).$ We have $P(e)=M(e)$ and
the functor $\theta_e$ is the identity functor. The composition or
the direct sum of two projective functors are again projective
functors. With respect to these two operations, projective
endofunctors on $\mc{O}_0$ are (up to isomorphism) generated by the
projective functors $\theta_i:=\theta_{s_i}$ corresponding to the
simple transpositions/reflections $s_i=(i,i+1)$. The functor
$\theta_i$ is called {\it the translation through the $i$-th wall}.
The functor $\theta_w$ is a direct summand of $\theta_{i_k}\dots
\theta_{i_1},$ for any reduced decomposition $w=s_{i_1}\dots
s_{i_k}.$ The induced endomorphism $[\theta_i]$ of the Grothendieck
group acts (in the basis given by Verma modules) by
$$[\theta_i(M(w))]=[M(w)]+[M(ws_i)].$$

Now we are prepared to explain the categorification. We first fix
the unique isomorphism $\varphi$ of groups
\begin{eqnarray*}
\varphi:\quad K(\mc{O}_0)&\longrightarrow&\Z[S_n]\\
\left[M(w)\right]&\longmapsto& w
\end{eqnarray*}
and define $C_w:=\varphi([P(w)])$. Then the $C_w$, $w\in W$, form a
basis $\mathbf{a}$ of $\mathbb{Z}[S_n]$.\\

The action of $[\theta_i]$ corresponds under $\varphi$ to the
endomorphism of $\Z[S_n]$ given by right multiplication with
$C_{s_i}:=1+s_i$.

The defining relations of the generators $1+s_i$ in
$\mathbb{Z}[S_n]$ lift to isomorphisms of functors as follows
\begin{eqnarray*}
\theta_i^2 & \cong & \theta_i \oplus \theta_i , \\
\theta_i\theta_j & \cong & \theta_j \theta_i \hspace{0.1in}
  \mathrm{if} \hspace{0.1in} |i-j|>1, \\
\theta_i\theta_{i+1}\theta_i \oplus \theta_{i+1} & \cong  &
\theta_{i+1}\theta_i \theta_{i+1}\oplus \theta_i.
\end{eqnarray*}
Here, the last isomorphism follows from the existence of
decompositions of functors
\begin{eqnarray*}
\theta_i\theta_{i+1}\theta_i & \cong & \theta_{w_1}  \oplus \theta_i, \\
\theta_{i+1}\theta_i \theta_{i+1} & \cong & \theta_{w_1}\oplus
\theta_{i+1},
\end{eqnarray*}
where $w_1=s_is_{i+1}s_i=s_{i+1}s_is_{i+1}.$ In particular,
$[\theta_{w_1}]$ corresponds under $\varphi$ to the right
multiplication with $1+s_i+s_{i+1}+s_is_{i+1}+s_{i+1}s_i
+s_is_{i+1}s_i.$

By the classification theorem of projective functors, the
endomorphism $[\theta_w]$, $w\in W$, corresponds then to right
multiplication with the element $C_w$. From this one can then
actually deduce that the multiplication in the basis $\mathbf{a}$
has non-negative integral coefficients
\begin{equation}\label{Sn}
 C_w C_{w'}=\sum_{w''} {c_{w,w'}^{w''}} C_{w''},
\hspace{0.2in} c_{w,w'}^{w''}\in \Z_{\geq0}.
\end{equation}
Hence we are in the situation of \eqref{cs} and are looking for an
abelian categorification of
$(\mathbb{Z}[S_n],\mathbf{a},\mathbb{Z}[S_n])$. We already have the
isomorphism $\varphi$ and the exact endofunctor $\theta_w$
corresponding to the generator $C_w$ satisfying condition
(C-\ref{I}).

The composition of two projective functors decomposed as a direct
sum of indecomposable functors $[\theta_w]$, $w\in W$, has
nonnegative integral coefficients, and the equations \eqref{Sn} turn
into isomorphisms of functors
\begin{equation}\label{equa1}
 \theta_w\theta_{w'}\cong\bigoplus_{w''}{(\theta_{w''})}^{c_{ww'}^{w''}},
\hspace{0.2in} c_{ww'}^{w''}\in \Z_{\geq0},
\end{equation}

It turns out that each $[\theta_w]$ acts by a multiplication with a
linear combination of $y$'s for $y\le w.$
Moreover, all coefficients are nonnegative integers.
For instance, if $w\in S_4$ then $[\theta_w]\doteq\sum_{y\le w} y$,
with two exceptions:
\begin{eqnarray*}
\begin{array}[thb]{ccll}
[\theta_w]&\doteq&\sum_{y \le w} y + 1 + s_2,
&w = s_2 s_1 s_3 s_2,\\
\left[\theta_w\right]&\doteq&\sum_{y\le w} y + 1 + s_1+ s_3
+s_1s_3 , &w = s_1s_3 s_2 s_1 s_3.
\end{array}
\end{eqnarray*}

We can summarise the above results into a theorem.

\begin{theorem} The action of the indecomposable projective functors
$\theta_w, w\in S_n$, on the block $\mc{O}_0$ for $\mf{sl}_n$
categorifies the right regular representation of the integral group
ring of the symmetric group $S_n$ (in the basis $\bf{a}$ of the
elements $C_w$, $w\in S_n$).
\end{theorem}

This theorem is due to Bernstein and Gelfand, see \cite{BeGe}, where it was stated in
different terms, since the word
``categorification'' was not in the mathematician's vocabulary back then.
In fact, Bernstein and Gelfand obtained a more general result  by considering
 any simple Lie algebra $\mf{g}$
instead of  $\mf{sl}_n$ and its Weyl group $W$ in place of $S_n.$

In the explanation to the theorem we did not give a very explicit
description of the basis $\bf{a}$ due to the fact that there is no
explicit (closed) formula for the elements $C_w$ available. However,
the elements $C_w$ can be obtained by induction (on the length of
$w$) using the Kazhdan-Lusztig theory \cite{KL1}, \cite{KL2}. The
Kazhdan-Lusztig theory explains precisely the complicated interplay
between the basis ${\bf a}$ and the standard basis of
$\Z[S_n]$.

\vspace{0.15in}


{\bf 3.}
Parabolic blocks of $\mc{O}$ categorify representations
of the symmetric group $S_n$ induced from the sign representation of
parabolic subgroups.

\vspace{0.15in}

Let $\mu=(\mu_1, \dots, \mu_k), \mu_1+\dots + \mu_k=n$, be a
composition of $n$ and $\lambda=(\lambda_1, \dots, \lambda_k)$ the
corresponding partition. In other words, $\lambda$ is a permutation
of the sequence $\mu$ with $\lambda_1\ge \lambda_2 \ge \dots \ge
\lambda_k.$ Denote by $p_{\mu}$ the subalgebra of $\mf{sl}_n$
consisting of $\mu$-block upper-triangular matrices. Consider the
full subcategory $\mc{O}^{\mu}$ of $\mc{O}$ which consists of all
modules $M$ on which the action of $U(p_{\mu})$ is locally finite.
The category $\mc{O}^{\mu}$ is an example of a \emph{parabolic}
subcategory of $\mc{O}$, introduced in \cite{Rocha-Caridi}. A simple
object $L(\lambda)$ of $\mc{O}$ belongs to $\mc{O}^{\mu}$ if and
only if the weight $\lambda$ is positive integral with respect to
all roots of the Lie algebra $p_{\mu}$. The two extreme cases are
$\mu=(1,1,\dots, 1),$ in which case $\mc{O}^{\mu}$ is all of
$\mc{O},$ and $\mu=(n),$ for $\mc{O}^{(n)}$ is the semisimple
category consisting exactly of all finite-dimensional
$\mf{sl}_n$-modules.

The direct sum decomposition (\ref{dir-eq}) induces a similar decomposition
of the  parabolic category:
$$ \mc{O}^{\mu} \cong \oplusop{\nu \in  \mf{h}^{\ast}/S_n} \mc{O}^{\mu}_{\nu}.$$
Each category $\mc{O}^{\mu}_{\nu}$ is either trivial (i.e. contains
only the zero module) or equivalent to the category of
finite-dimensional modules over some finite-dimensional $\C$-algebra
(but describing this algebra explicitly for interesting $\mu$ and
$\nu$ is a hard problem, see \cite{BoeN}). Unless $\mu=(1^n),$ for
generic $\nu$ the summand $\mc{O}^{\mu}_{\nu}$ is trivial. Again,
the most complicated summands are the $\mc{O}^{\mu}_{\nu}$ where the
orbit $\nu$ contains a dominant regular integral weight. Translation
functors establish equivalences between such summands for various
such $\nu,$ and allow us to restrict our consideration to the block
$\mc{O}_{0}^{\mu}$ corresponding to the (shifted) orbit through $0$.
The inclusion
$$\mc{O}_{0}^{\mu}\subset \mc{O}_0$$
is an exact functor and induces an
inclusion of Grothendieck groups
\begin{equation}\label{inc-gr}
K(\mc{O}_{0}^{\mu}) \subset K(\mc{O}_0).
\end{equation}
Indeed, the Grothendieck group of $\mc{O}_0$ is free abelian with
generators $[L(w)],$ $w\in S_n.$ A simple module $L(w)$ lies in
$\mc{O}_0^{\mu}$ if and only if $w$ is a minimal left coset
representative for the subgroup $S_{\mu}$ of $S_n$ (we informally
write $w\in (S_{\mu}\backslash S_n)_{short}$). The Grothendieck
group of $\mc{O}_0^{\mu}$ is then the subgroup of $K(\mc{O}_0)$
generated by such $L(w).$

The analogues of the Verma modules in the parabolic case are the so-called
parabolic Verma modules
$$M(p_{\mu},V) = U(\mf{sl}_n) \otimes_{U(p_{\mu})}V,$$
where $V$ is a finite-dimensional simple $p_{\mu}$-module. The
module $M(p_{\mu},V)$ is a homomorphic image of some ordinary Verma
module from $\mc{O}$, in particular, it has a unique simple quotient
isomorphic to some $L(w)$ for some unique $w\in S_n$. In this way we
get a canonical bijection between parabolic Verma modules in
$\mc{O}^{\mu}_0$ and the set $(S_{\mu}\backslash S_n)_{short}$ of
shortest coset representatives. Hence it is convenient to denote the
parabolic Verma module with simple quotient $L(w)$, $w\in
(S_{\mu}\backslash S_n)_{short}$, simply by $M^{\mu}(w).$

Generalised Verma modules provide a basis for the Grothendieck group
of $\mc{O}^{\mu}_0.$ Under the inclusion (\ref{inc-gr}) of
Grothendieck groups the image of the generalised  Verma module
$M^{\mu}(w)$ is the alternating sum of Verma modules, see
\cite{Rocha-Caridi} and \cite{Lepow}:
\begin{equation} \label{alt-verm}
 [M^{\mu}(w)] = \sum_{u\in S_{\mu}} (-1)^{l(u)} [M(uw)].
\end{equation}

Since the projective endofunctors $\theta_w$ preserve
$\mc{O}^{\mu}_0,$ the inclusion (\ref{inc-gr}) is actually an
inclusion of $S_n$-modules, and, in view of the formula
(\ref{alt-verm}), we can identify $K(\mc{O}^{\mu}_0)$ with the
submodule $I^-_{\mu}$ of the regular representation of $S_n$
isomorphic to the representation induced from the sign
representation of $S_{\mu},$
$$ I^-_{\mu} \cong \mathrm{Ind}_{S_{\mu}}^{S_n}  \Z v ,$$
where we denoted by $\Z v$ the sign representation, so that $w v =
(-1)^{l(w)} v$ for $w\in S_{\mu}.$

To summarise, we have:

\begin{theorem}
\label{perm} The action of the projective functors $\theta_w$, $w\in
W$, on the parabolic subcategory $\mc{O}_0^{\mu}$ of $\mc{O}$
categorifies the induced representation $I_{\mu}^-$ of the integral
group ring of the symmetric group $S_n$ (with basis ${\bf
a}=\{C_w\}_{w\in S_n}$).
\end{theorem}

As in the previous example the Grothendieck group
$K(\mc{O}^{\mu}_0)$ has three distinguished basis, given by simple
objects, projective objects, and parabolic Verma modules
respectively.

\vspace{0.2in}

\emph{Remark:} If we choose a pair $\mu, \mu'$ of decompositions
giving rise to the same partition $\lambda$ of $n$, then the modules
$I_{\mu}^-$ and $I_{\mu'}^-$ are isomorphic, and will be denoted
$I_{\lambda}^-.$ However, the categories $\mc{O}^{\mu}$ and
$\mc{O}^{\mu'}$ are not equivalent in general, which means the two
categorifications of the induced representation $I_{\lambda}^-$ are
also not equivalent. This problem disappears if we leave the world
of abelian categorifications, since the derived categories
$D^b(\mc{O}^{\mu})$ and $D^b(\mc{O}^{\mu'})$ are equivalent
\cite{Kh-springer}. The equivalence is based on the geometric
description of $\mc{O}^{\mu}$ and $\mc{O}^{\mu'}$ in terms of
complexes of sheaves on partial flag varieties.

\vspace{0.2in}


{\bf 4.} Self-dual projectives in a parabolic block categorify irreducible
representations of the symmetric group.

\vspace{0.1in}

Let $I_{\mu}$ be the representation of $\Z[S_n]$ induced from the
trivial representation of the subgroup $S_{\mu}.$ Up to isomorphism,
it only depends on the partition $\lambda$ associated with $\mu$.
Partitions of $n$
naturally index the isomorphism classes of irreducible
representations of $S_n$ over any field of characteristic zero (we
use $\Q$ here). Denote by $\mathcal{S}_{\Q}(\lambda)$ the
irreducible (Specht) module associated with $\lambda.$ It is an
irreducible representation defined as the unique common irreducible
summand of $I_{\lambda}\otimes \Q$ and $I_{\lambda^{\ast}}^-\otimes
\Q,$ where $\lambda^{\ast}$ is the dual partition of $\lambda.$
Passing to duals, we see that $\mathcal{S}_{\Q}(\lambda^{\ast})$ is
the unique common irreducible summand of $I_{\lambda^{\ast}}\otimes
\Q$ and $I_{\lambda}^-\otimes \Q.$

We have already categorified the representation $I_{\lambda}^-$ (in
several ways) via the parabolic categories $\mathcal{O}_0^{\mu},$
where $\mu$ is any decomposition for $\lambda.$ It's natural to try
to realise a categorification of some integral lift
$\mathcal{S}(\lambda^{\ast})$ of the irreducible representation
$\mathcal{S}_{\Q}(\lambda^{\ast})$ via a suitable subcategory of
some  $\mathcal{O}_0^{\mu}$ stable under the action of projective
endofunctors.

The correct answer, presented in \cite{KMS}, is to pass to a
subcategory generated by those projective objects in
$\mathcal{O}_0^{\mu}$ which are also injective. Note that these
modules are neither projective nor injective in $\cO_0$ (unless if
$\cO_o^\mu=\cO_0).$

Any projective object in  $\mathcal{O}_0^{\mu}$ is isomorphic to a
direct sum of indecomposable projective modules $P^{\mu}(w),$ for
$w\in (S_{\mu}\backslash S_n)_{short}.$ Let $J\subset
(S_{\mu}\backslash S_n)_{short}$ be the subset indexing
indecomposable projectives modules that are also injective: $w\in J$
if and only if $P^{\mu}(w)$ is injective. Projective endofunctors
$\theta_w,$ $w\in S_n$, take projectives to projectives and
injectives to injectives. Therefore, they take projective-injective
modules (modules that are both projective and injective, also called
self-dual projective, for instance, in Irving \cite{Ir1}) to
projective-injective modules.

The category of projective-injective modules is additive, not abelian.
To remedy this, consider the full subcategory $\mc{C}^{\mu}$
of $\mc{O}_0^{\mu}$ consisting of modules $M$ admitting
a resolution
\begin{equation}\label{formula1}
 P_1\lra P_0 \lra M \lra 0
\end{equation}
with projective-injective $P_1$ and $P_0.$ The category
$\mc{C}^{\mu}$ is abelian and stable under all endofunctors
$\theta_w$ for $w\in S_n$, see \cite{KMS}.

Irving \cite{Ir1} classified projective-injective modules in
$\mc{O}_0^{\mu}.$ His results were interpreted in \cite{KMS}
in the language of categorification:

\begin{theorem}
\label{Specht} The action of the projective endofunctors $\theta_w$,
$w\in S_n$, on the abelian category $\mc{C}^{\mu}$ categorifies (after
tensoring the Grothendieck group with $\Q$ over $\Z$)
the irreducible representation
$\mc{S}_{\Q}(\lambda^{\ast})$  of the symmetric group $S_n.$
\end{theorem}

\vspace{0.1in}

The Grothendieck group $K(\mc{C}^{\mu})$ is a module
over the integral group ring of $S_n,$ with $s_i$ acting by $[\theta_i]-\mathrm{Id},$
and the theorem says that $K(\mc{C}^{\mu})\otimes_{\Z} \Q $
is an irreducible representation of the symmetric group corresponding to
the partition $\lambda^{\ast}.$
Several explicit examples of categorifications via $\mc{C}^{\mu}$
will be given in Section~\ref{four-examples}.

\vspace{0.1in}

\emph{Remark:} Suppose $\mu$ and $\nu$ are two decompositions of the
same partition $\lambda.$ It's shown in \cite{MS-projinj} (Theorem
5.4.(2)) that the categories $\mc{C}^{\mu}$ and $\mc{C}^{\nu}$ are
equivalent, through an equivalence which commutes with the action of
the projective functors $\theta_u$ on these categories (the
equivalence is given by a non-trivial composition of derived
Zuckerman functors). Therefore, the categorification of
$\mc{S}(\lambda^{\ast})$ does not depend on the choice of the
decomposition $\mu$ that represents $\lambda,$ and we can denote the
category $\mc{C}^{\mu}$ by $\mc{C}^{\lambda}.$ (This should be
compared with the remark after Theorem~\ref{perm}.)

\vspace{0.1in}

\emph{Remark:} Theorem~\ref{perm} and Theorem~\ref{Specht} can be
generalised to arbitrary semi-simple complex finite-dimensional Lie
algebras, see \cite{MSnew}. However, in the general case
Theorem~\ref{Specht} does not categorify simple modules for the
corresponding Weyl group but rather the Kazhdan-Lusztig cell modules
from \cite{KL1}. This can be used to describe the so-called
``rough'' structure of generalised Verma modules, which shows that
``categorification theoretic'' ideas can lead to new
results in representation theory.

\vspace{0.1in}

\emph{Remark:} The inclusion of categories $\mc{C}^{\mu}\subset
\mc{O}^{\mu}_0$ is not an exact functor; however, it is a part of a
very natural filtration of the category $\mc{O}^{\mu}_0$
which can be defined using
the Gelfand-Kirillov dimension of modules, see \cite[6.9]{MSnew}. To
get the inclusion of Grothendieck groups analogous to the inclusion
of representations from the irreducible Specht module into the
induced sign representation, we pass to the subgroup
$K'(\mc{C}^{\mu})$ of $K(\mc{C}^{\mu})$ generated by the images of
projective modules in $\mc{C}^{\mu}.$ This additional technicality
is necessary as the category $\mc{C}^{\mu}$ does not have finite
homological dimension in general. The subgroup $K'(\mc{C}^{\mu})$ is
 always a finite index subgroup, stable under the action of the $[\theta_w]$'s.
We denote this subgroup by $\mc{S}'(\lambda^{\ast}):$
$$\mc{S}'(\lambda^{\ast}) \define K'(\mc{C}^{\mu}) \subset
K(\mc{C}^{\mu}) \cong  \mc{S} (\lambda^{\ast}).$$
The inclusion of categories $\mc{C}^{\mu}\subset \mc{O}^{\mu}_0$
induces the inclusion
$$\mc{S}'(\lambda^{\ast}) \subset K(\mc{O}^{\mu}_0) \cong
 I_{\mu}^-$$ of $\Z[S_n]$-modules, hence realising the integral lift $\mc{S}'(\lambda^{\ast})$
of the Specht module as a subrepresentation of $I_{\mu}^-.$

\vspace{0.2in}


{\bf 5.} Categorification of the induced representations $I_{\mu}$
via projectively presentable modules.

\vspace{0.2in}

Let $\mathcal{P}^{\mu}$ denote the category of all modules $M$
admitting a resolution \eqref{formula1} in which each indecomposable
direct summand of both $P_0$ and $P_1$ has the form $P(w)$, where
$w$ is a longest  left coset representative for $S_{\mu}$ in $S_n$
(we will write $w\in (S_{\mu}\backslash S_n)_{long}$). Such modules
are called $p_{\mu}$-presentable modules, see \cite{MS-shuffling}.
As in the previous example, the category $\mathcal{P}^{\mu}$ is
stable under all endofunctors $\theta_w$, $w\in S_n$.

By definition, $\mathcal{P}^{\mu}$ is a subcategory of $\mc{O}_0$,
but just as in the example above, the natural inclusion functor is
not exact. The category $\mc{P}^{\mu}$ does not have finite
homological dimension in general, so we again pass to the subgroup
$K'(\mathcal{P}^{\mu})$  of $K(\mathcal{P}^{\mu})$, generated by the
images of indecomposable projective modules in $\mathcal{P}^{\mu}$.
The latter are (up to isomorphism) the  $P(u)$,  $u\in
(S_{\mu}\backslash S_n)_{long}$. This is a finite index subgroup,
stable under the action of the $[\theta_w]$'s and we have the
following statement proved in \cite{MS-shuffling}:

\begin{theorem} The
action of projective endofunctors on the abelian category
$\mc{P}^{\mu}$ categorifies (after tensoring with $\Q$ over $\Z$)
the induced representation $(I_{\mu})_{\Q}$ of the group
algebra of the symmetric group $S_n$ (with the basis ${\bf
a}=\{C_w\}_{w\in S_n})$.
\end{theorem}

\vspace{0.1in}

Consider the diagram of $\Q[S_n]$-modules
$$(I_{\lambda})_{\Q}\stackrel{\iota_1}{\longrightarrow} \Q[S_n]
\stackrel{p_1}{\longrightarrow} (I_{\lambda^{\ast}}^-)_{\Q} .$$ The
map $\iota_1$ is the symmetrisation inclusion map, while $p_1$ is
the antisymmetrization quotient map. We have
$$\mathcal{S}_{\Q}(\lambda)\define p_1 \iota_1 ((I_{\lambda})_{\Q}).$$
The map $\iota_1$ is categorified as the inclusion of
$\mathcal{P}^{\mu}$ to $\mc{O}_0$. The map $p_1$ is categorified as
the projection of $\mc{O}_0$ onto $ \mc{O}_0^{\mu^*}$, where $\mu^*$
is some decomposition corresponding to $\lambda^*$. Unfortunately,
the composition of the two functors categorifying these two maps
will be trivial in general. To repair the situation we first project
$\mc{P}^{\mu}$ onto the full subcategory of $\mc{P}^{\mu}$ given by
simple objects of minimal possible Gelfand-Kirillov dimension. It is
easy to see that the image category contains enough projective
modules, and using the equivalence constructed in
\cite[Theorem~I]{MSnew}, these projective modules can be
functorially mapped to projective modules in $\mathcal{C}^{\mu^*}$,
where $\mu^*$ is a (good choice of) composition with associated
partition $\lambda^*$. The latter category embeds into
$\mc{O}_0^{\mu^*}$ as was explained in the previous example.


\vspace{0.2in}

{\bf 6.} The representation theory of groups like
$GL(n,\mathbb{C})$, considered as a real Lie group, naturally leads
to the notion of Harish-Chandra bimodules. A Harish-Chandra bimodule
over $\mf{sl}_n$ is a finitely-generated module over the universal
enveloping algebra $U(\mf{sl}_n\times \mf{sl}_n)$ which decomposes
into a direct sum of finite-dimensional $U(\mf{sl}_n)$-modules with
respect to the diagonal copy $\{(X,-X)| X\in \mf{sl}_n\}$ of
$\mf{sl}_n.$ Let $\mc{HC}_{0,0}$ be the category of Harish-Chandra
bimodules which are annihilated, on both sides, by some power of the
maximal ideal $I_0$ of the centre $Z$ of $U(\mf{sl}_n).$ Here $I_0$
is the annihilator of the trivial $U(\mf{sl}_n)$-module considered
as a $Z$-module. Thus, $M\in \mc{HC}_{0,0} $ if and only if $ xM = 0
= Mx$ for all $x\in I_0^N$ for $N$ large enough.

By \cite[Section 5]{BeGe}  there exists an exact and fully faithful functor
$$\mc{O}_0\lra \mc{HC}_{0,0}.$$
Moreover, this functor induces an isomorphism of Grothendieck groups
 $$ K(\mc{O}_0) \cong K(\mc{HC}_{0,0}).$$
Since the former group is isomorphic to $\Z[S_n],$ we can identify
the Gro\-then\-di\-eck group of $\mc{HC}_{0,0}$ with $\Z[S_n]$ as well.

The advantage of bimodules is that we now have two sides and can
tensor with a finite-dimensional $\mf{sl}_n$-module both on the left
and  on the right. In either case, we preserve the category of
Harish-Chandra bimodules. Taking all possible direct summands of
these functors and restricting to endofunctors on the subcategory
$\mc{HC}_{0,0}$ leads to two sets of commuting projective functors,
$\{\theta_{r,w}\}_{w\in S_n}$ and $\{\theta_{l,w}\}_{w\in S_n}$
which induce endomorphisms on the Grothendieck group
$K(\mc{HC}_{0,0})\cong \Z[S_n]$ given by left and right
multiplication with $\{C_w\}_{w\in S_n}$ respectively. Summarising,
we have

\begin{theorem} The action of the functors $\{\theta_{r,w}\}_{w\in S_n}$ and $\{\theta_{l,w}\}_{w\in S_n}$ on the category $\mc{HC}_{0,0}$ of Harish-Chandra
bimodules for $\mf{sl}_n$ with generalised trivial character on both
sides categorifies $\Z[S_n],$ viewed as a bimodule over itself. The
functors $\{\theta_{r,w}\}_{w\in S_n}$ induces the left
multiplication with $C_w$ on the Grothendieck group, whereas the
functors $\{\theta_{l,w}\}_{w\in S_n}$ induces the right
multiplication with $C_w$ on the Grothendieck group,
\end{theorem}

The first half of the theorem follows at once from \cite{BeGe}, the
second half from \cite{Soergel-HC}. The category of Harish-Chandra
bimodules is more complicated than the category $\mc{O}$. For
example, $\mc{HC}_{0,0}$ does not have enough projectives, and is not Koszul
with respect to the natural grading, in contrast to $\cO_0$ (for the Koszulity of
$\cO_0$ see \cite{BGS}). The study of translation functors on Harish-Chandra
modules goes back to Zuckerman~\cite{Zuckerman}.

\vspace{0.2in}


{\bf 7.} In the following we will mention several instances of
categorifications of modules over Lie algebras. Our definition of
categorification required an associative algebra rather than a Lie
algebra, so one should think of this construction as a
categorification of representations of the associated universal
enveloping algebra.

Let  $V$ be the fundamental two-dimensional representation of the
complex Lie algebra $\mf{sl}_2.$ Denote by $\{e,f,h\}$ the standard
basis of $\mf{sl}_2.$ The $n$-th tensor power of $V$ decomposes into
a direct sum of weight spaces:
$$ V^{\otimes n} = \oplusoop{k=0}{n} V^{\otimes n}(k) ,$$
where $hx= (2k-n)x$ for $x\in V^{\otimes n}(k).$

A categorification of $V^{\otimes n}$ was constructed in \cite{BFK}.
The authors considered certain singular blocks $\mc{O}_{k,n-k}$ of
the category $\mc{O}$ for $\mf{sl}_n.$ The Grothendieck group of
this block has rank $\binom{n}{k}$ equal to the dimension of the
weight space $V^{\otimes n}(k) ,$ and there are natural isomorphisms
$$ K(\mc{O}_{k,n-k}) \otimes_{\Z} \C \cong V^{\otimes n} (k) .$$
The Grothendieck group of the direct sum
$$ \mc{O}_n = \oplusoop{k=0}{n}\mc{O}_{k,n-k}$$
is isomorphic to $V^{\otimes n}$ (after tensoring with $\C$ over $\Z$).
Suitable translation functors $\mc{E}, \mc{F}$ in $\mc{O}_n$
lift the action of the generators $e,f$ of $\mf{sl}_2$ on $V^{\otimes n}.$

To make this construction compatible with Definition~\ref{def-1} one
should switch to Lusztig's version $\stackrel{\cdot}{\mathbf{U}}$ of
the universal enveloping algebra $U(\mf{sl}_2)$ (see \cite{L},
\cite{BFK}) and set $q=1$. Instead of the unit element, the ring
$\stackrel{\cdot}{\mathbf{U}}$ contains idempotents  $1_n, n\in \Z,$
which can be viewed as projectors onto integral weights. The Lusztig
basis $\stackrel{\cdot}{\mathbb{B}}$ in
$\stackrel{\cdot}{\mathbf{U}}$ has the positivity property required
by Definition~\ref{def-1}, and comes along with an integral version
$V^{\otimes n}_{\Z}$ of the tensor power representation. The triple
$(\stackrel{\cdot}{\mathbf{U}}, \stackrel{\cdot}{\mathbb{B}},
V^{\otimes n}_{\Z})$ is categorified using the above-mentioned
category $\mc{O}_n$ and projective endofunctors of it. In fact, each
element of $\stackrel{\cdot}{\mathbb{B}}$ either corresponds to an
indecomposable projective endofunctor on $\mc{O}_n$ or acts by $0$
on $V^{\otimes n}_{\Z}.$ We refer the reader to \cite{BFK} for
details, to \cite{ChRo} for an axiomatic development of $\mf{sl}_2$
categorifications, and to \cite{BFK} and \cite{Catharina-categ} for
a categorification of the Temperley-Lieb algebra action on
$V^{\otimes n}_{\Z}$ via projective endofunctors on the category
Koszul dual to $\mc{O}_n$ (see \cite{BGS} and \cite{MOS} for details
on Koszul duality).

\vspace{0.06in}

The Lie algebra $\mf{sl}_2$ has one irreducible $(n+1)$-dimensional
representation $V_n,$ for each $n\ge 0$ ($V_1\cong V,$ of course). A
categorification of arbitrary tensor products $V_{n_1}\otimes \dots
\otimes V_{n_m}$ is described in \cite{FKS}. This tensor product is
a submodule of $V^{\otimes n},$ where $n=n_1+\dots +n_m.$ Knowing
that $\mc{O}_n$ categorifies $V^{\otimes n},$ we find a
"subcategorification," a subcategory of $\mc{O}_n$ stable under the
action of projective functors, with the Grothendieck group naturally
isomorphic to  $V_{n_1}\otimes \dots \otimes V_{n_m}.$ The
subcategory has an intrinsic description via Harish-Chandra modules
similar to the one from Example {\bf 5}.

\vspace{0.15in}


{\bf 8.} A categorification of arbitrary tensor products of
fundamental representations $\Lambda^iV,$ where $V$ is the
$k$-dimensional $\mf{sl}_k$-representation  and $1\le i\le k-1$ was
found by J.~Sussan \cite{Sussan}. A tensor product
$\Lambda^{i_1}V\otimes \dots \otimes\Lambda^{i_r}V$ decomposes into
weight spaces $\Lambda^{i_1}V\otimes \dots
\otimes\Lambda^{i_r}V(\nu),$ over various integral weights $\nu$ of
$\mf{sl}_k.$ Each weight space becomes the Grothendieck group of a
parabolic-singular block of the highest weight category for
$\mf{sl}_N,$ where $N=i_1+\dots+i_r.$ For the parabolic subalgebra
one takes the Lie algebra of traceless $N\times N$ matrices which
are $(i_1, \dots, i_r)$ block upper-triangular. The choice of the
singular block is determined by $\nu.$ Translation functors between
singular blocks, restricted to the parabolic category, provide an
action of the generators $\mc{E}_j$ and $\mc{F}_j$ of the Lie
algebra $\mf{sl}_k.$ Relations in the universal enveloping algebra
lift to functor isomorphisms. Conjecturally, Sussan's
categorification satisfies the framework of Definition~\ref{def-1}
above, with respect to Lusztig's completion
$\stackrel{\cdot}{\mathbf{U}}$ of the universal enveloping algebra
of $\mf{sl}_k$ and Lusztig's canonical basis there.

\vspace{0.15in}


{\bf 9.} Ariki, in a remarkable paper \cite{Ariki-dec}, categorified
all finite-dimensional irreducible representations of $\mf{sl}_m,$
for all $m,$ as well as integrable irreducible representations of
affine Lie algebras $\widehat{\mf{sl}}_r.$ Ariki considered certain
finite-dimensional quotient algebras of the affine Hecke algebra
$\widehat{H}_{n,q},$ known as Ariki-Koike cyclotomic Hecke algebras,
which depend on a number of discrete parameters. He identified the
Grothendieck groups of blocks of these algebras, for generic values
of $q\in \C,$ with the weight spaces $V_{\lambda}(\mu)$ of
finite-dimensional irreducible representations
$$V_{\lambda} =\oplusop{\mu} V_{\lambda}(\mu)$$
of $\mf{sl}_m.$ Direct summands
of the induction and restriction functors between cyclotomic Hecke algebras
for $n$ and $n+1$ act on the Grothendieck group as generators $e_i$ and
$f_i$ of $\mf{sl}_m.$

Specialising $q$ to a primitive $r$-th root of unity, Ariki obtained
a categorification of integrable irreducible representations of the
affine Lie algebra $\widehat{\mf{sl}}_r.$

We conjecture that direct summands of arbitrary compositions of
Ariki's induction and restriction functors correspond to elements of
the Lusztig canonical basis $\stackrel{\cdot}{\mathbb{B}}$ of
Lusztig's completions $\stackrel{\cdot}{\mathbf{U}}$ of these
universal enveloping algebras. This conjecture would imply that
Ariki's categorifications satisfy the conditions of Definition 1.

Lascoux, Leclerc and Thibon, in an earlier paper \cite{LLT},
categorified level-one irreducible
$\widehat{\mf{sl}}_r$-representations, by identifying them with the
direct sum of Grothendieck groups of finite-dimensional Hecke
algebras $H_{n,q},$ over all $n\ge 0,$ with $q$ a primitive
 $r$-th root of unity. Their construction is a special case of Ariki's.
We also refer the reader to related works \cite{Ariki-book},
\cite{Grojnowski-slp}. Categorifications of the adjoint
representation and of irreducible $\mf{sl}_m$-representations with
highest weight $\omega_j+\omega_k$ are described explicitly in \cite{HK1},
\cite{HK2} and \cite{Chen}.

Another way to categorify all irreducible finite-dimensional representations
of $\mf{sl}_m$, for all $m,$ was found by Brundan
and Kleshchev~\cite{BK1}, via the representation theory of W-algebras.
There is a good chance that their categorification is equivalent to that of Ariki,
and that an equivalence of two categorifications can be constructed
along the lines of Arakawa-Suzuki~\cite{AS} and Brundan-Kleshchev~\cite{BK2}.

\vspace{0.2in}


{\bf Biadjointness.}
Definition 1 of (weak) categorifications was minimalistic. Categorifications
in the above examples share extra properties, the most prominent of which
is biadjointness: there exists an involution $a_i\to a_{i'}$ on the basis $\ba$ of $A$
such that the functor $F_{i'}$ is both left and right adjoint to $F_i.$ This
is the case in the examples {\bf 2} through {\bf 9},
while in example {\bf 1} the functors are almost biadjoint. Namely,
the induction functor $F_x$ lifting the action of $x$ is left adjoint to
the restriction functor $F_{\partial}$ (which lifts the action of $\partial$)
and right adjoint to $F_{\partial}$ conjugated by an involution.

A conceptual explanation for the pervasiveness of biadjointness in
categorifications is given by the presence of the $\mathrm{Hom}$
bifunctor in any abelian category. The $\mathrm{Hom}$ bifunctor in
$\mc{B}$ descends to a bilinear form on the Grothendieck group $B$
of $\mc{B},$ via
$$([M],[N]) \define \dim \mathrm{Hom}_{\mc{B}}(M,N),$$
where  $M$ is projective or $N$ is injective, and some standard technical
conditions are satisfied. When a representation naturally
comes with a bilinear form, the form is usually compatible
with the action of $A$: there exists an involution $a \to a'$ on $A$
such that $(ax,y)=(x,a'y)$ for $x,y\in B.$ A categorification of this
equality should be an isomorphism
$$ \mathrm{Hom}(F_a M,N) \cong \mathrm{Hom}(M,F_{a'}N)$$
saying that the functor lifting the action of $a'$ is right adjoint to
the functor lifting the action of $a.$ If the bilinear form is symmetric,
we should have the adjointness property in the other direction as well,
leading to biadjointness of $F_a$ and $F_{a'}.$

\vspace{0.06in}

A beautiful approach to $\mf{sl}_2$ categorifications via
biadjointness was developed by Chuang and Rouquier \cite{ChRo} (see also 
\cite{Rouquier-survey}). The
role of biadjointness in TQFTs and their categorifications is
clarified in \cite[Section 6.3]{Kh-FVIT}. An example how the
existence of a categorification with a bilinear form can be used to
determine dimensions of hom-spaces can be found in
\cite{Catharina-corners}.

\vspace{0.2in}


{\bf Grading and $q$-deformation.} In all of the above examples, the
data $(A,\ba, B)$ that is being categorified admits a natural
$q$-deformation $(A_q, \ba_q, B_q).$ Here $A_q$  is a
$\Z[q,q^{-1}]$-algebra, $B_q$ an $A_q$-module, and $\ba_q$ a basis
of $A_q.$ We assume that both $A_q$ and $B_q$ are free
$\Z[q,q^{-1}]$-modules, that the multiplication in $A_q$ in the
basis $\ba_q$ has all coefficients in $\N[q,q^{-1}],$ and that
taking the quotient by the ideal $(q-1)$ brings us back to the
original data:
$$ A= A_q/(q-1)A_q, \hspace{0.2in} B= B_q/(q-1)B_q,
\hspace{0.2in} \ba_q \stackrel{q=1}{\lra} \ba.$$ An automorphism
$\tau$ of an abelian category $\mc{B}$ (more accurately, an
invertible endofunctor on $\mc{B}$) induces a $\Z[q,q^{-1}]$-module
structure on the Grothendieck group $K(\mc{B}).$ Multiplication by
$q$ corresponds to the action of $\tau$:
 $$ [\tau(M)] = q[M], \hspace{0.2in} [\tau^{-1}(M)] = q^{-1}[M].$$
In many of the examples, $\mc{B}$ will be the category of graded
modules over a graded algebra, and $\tau$ is just the functor which
shifts the grading. To emphasize this, we denote $\tau$ by $\{1\}$
and its $n$-th power by $\{n\}.$

\begin{definition} A (weak) abelian categorification of $(A_q,\ba_q,B_q)$ consists
of an abelian category $\mc{B}$ equipped with an invertible
endofunctor $\{1\},$  an isomorphism of $\Z[q,q^{-1}]$-modules
$\varphi:K(\mc{B})\overset{\sim}{\longrightarrow} B_q$ and exact
endofunctors $F_i: \mc{B} \lra \mc{B}$ that commute with $\{1\}$ and
such that the following hold
\begin{enumerate}[(qC-I)]
\item
$F_i$ lifts the action of $a_i$ on the module $B_q$, i.e. the
action of $[F_i]$ on the Grothendieck group corresponds to the
action of $a_i$ on $B_q,$ under the isomorphism $\phi$, in the sense
that the diagram below is commutative.
 $$\begin{CD}
    K(\mc{B})  @>{[F_i]}>>     K(\mc{B})   \\
   @V{\varphi}VV                   @VV{\varphi}V   \\
    B_q   @>{a_i}>>     B_q
  \end{CD} $$
\item  There are isomorphisms of functors
$$ F_i F_j \cong \bigoplus_{k} {F_k}^{c_{ij}^k},$$
i.e., the composition $F_i F_j$  decomposes as the direct sum of functors
$F_k$ with multiplicities $c_{ij}^k\in \N[q,q^{-1}]$
\end{enumerate}
\end{definition}

\vspace{0.1in}

The graded versions are well-known in all of the examples above
up to Example {\bf 8}.  In Example {\bf 1} the nilCoxeter algebra
$R_n$ is naturally
graded with $\deg(Y_i)=1.$ The inclusion $R_n\subset R_{n+1}$ induces
induction and restriction functors between categories of graded $R_n$
and $R_{n+1}$-modules. In the graded case, induction and restriction
functors satisfy the isomorphism
$$DX \cong XD\{1\}\oplus \mathrm{Id}$$
which lifts the defining relation $\partial x = q x\partial +1$ of
the $q$-Weyl algebra (see \cite{Kh-nilcoxeter} for more detail).

\vspace{0.1in}

An accurate framework for graded versions of examples {\bf 2}--{\bf
8} is a rather complicated affair. To construct a canonical grading
on a regular block of the highest weight category \cite{BGS}
requires \'etale cohomology, perverse sheaves \cite{BBD}, and the
Beilinson-Bernstein-Brylinski-Kashiwara localisation theorem
\cite{BB}, \cite{BK}. Soergel's approach to this grading is more
elementary \cite{Soergel-garben}, \cite{Soergel-gr},
\cite{Soergel-HC}, but still relies on these hard results. Extra
work is needed to show that translation or projective functors can
be lifted to endofunctors in the graded category \cite{Catharina-O}.

\vspace{0.1in}

Ariki's categorification of irreducible integrable representations
(Example {\bf 9} above)  should admit a graded version as well.

\vspace{0.2in}


\section{Four examples of categorifications of irreducible representations}
\label{four-examples}

In the example {\bf 4} above we categorified an integral lift of the
 irreducible representation $\mc{S}_{\Q}(\lambda^{\ast})$ of the symmetric group
via the abelian category $\mc{C}^{\lambda}$ built out of projective-injective
modules in a parabolic block of $\mc{O}.$
The category $\mc{C}^{\lambda}$ is equivalent to the category
of finite-dimensional representations over a finite-dimensional
algebra $A^{\lambda},$ the algebra of endomorphisms of the direct
sum of indecomposable projective-injective modules $P^{\mu}(w).$
Under this equivalence, projective
functors $\theta_i$ turn into the functors of tensoring with
certain $A^{\lambda}$-bimodules. It's not known how to
describe $A^{\lambda}$ and these bimodules explicitly, except in a few
cases, four of which are discussed below.

\vspace{0.15in}


{\bf a. The sign representation.} The sign representation of the
symmetric group (over $\Z$) is a free abelian group $\Z v$ on one
generator $v,$ with $s_i v= -v$ for all $i.$ It corresponds to the
partition $(1^n)$ of $n,$ which in our notation is $\lambda^{\ast}$
for $\lambda=(n).$ The parabolic category $\mc{O}_0^{(n)}$ has as
objects exactly the finite-dimensional modules from $\mc{O}_0$ since
the parabolic subalgebra in this case is all of $\mf{sl}_n.$

Actually, $\mc{O}_0$ has only one simple module with this property, the one-dimensional
trivial representation $\C.$ In our notation, this is
the module $L(e),$ the simple quotient of the Verma module $M(e)$ assigned to
the unit element of the symmetric group.

Consequently, any object of $\mc{O}_0^{(n)}$ is isomorphic to
a direct sum of copies of $L(e),$ and the category is semisimple.
Furthermore, the category $\mc{C}^{(n)}$ is all of $\mc{O}_0^{(n)}.$
Thus, $\mc{C}^{(n)}$ is equivalent to the category of finite-dimensional
$\C$-vector spaces. Projective functors $\theta_w$ act by zero
on $\mc{C}^{(n)}$ for all $w\in S_n, w\not= e,$ while $\theta_e$ is
the identity functor.

The graded version $\mc{C}_{gr}^{(n)}$ is equivalent to the
category of graded finite-dimensional $\C$-vector spaces. Again,
projective functors $\theta_w, w\not=e,$ act by zero, and $\theta_e$
is the identity functor.

Thus,  our categorification of the sign representation is rather trivial.

\vspace{0.2in}


{\bf b. The trivial representation.}
The trivial representation $\Z z$ of $\Z[S_n]$ is a free abelian
group on one generator $z,$ with the action $w z =z,$ $w\in S_n.$
The corresponding partition is $(n),$ with the dual partition
$\lambda=(1^n).$ Only one decomposition $\mu=(1^n)$
corresponds to the dual partition;  the parabolic subalgebra associated
with $(1^n)$ is the positive Borel
subalgebra, and the parabolic category $\mc{O}^{(1^n)}_0$ is all of $\mc{O}_0.$

The unique self-dual indecomposable projective $P$ in $\mc{O}_0$ is
usually called the big projective module. Its endomorphism algebra
$\mathrm{End}_{\mc{O}}(P)$ is isomorphic to the cohomology ring
$\mathrm{H}_n$ of the full flag variety $\operatorname{Fl}$ of
$\C^n,$ see \cite{Soergel-garben}.

The category $\mc{C}^{(1^n)}$ is equivalent to the category
of finite-dimensional $\mathrm{H}_n$-modules. The unique (up to isomorphism)
simple $\mathrm{H}_n$-module generates the Grothendieck group
$K_0(\mathrm{H-mod})\cong \Z.$

To describe how the functors $\theta_i$ act on $\mc{C}^{\lambda}$
consider generalised flag varieties
\begin{eqnarray*}
 \operatorname{Fl}_i & = & \{ 0\subset L_1\subset L_2\subset \dots \subset L_{n-1}\subset \C^n,
L_i' |   \\
& &
\dim(L_j)=j, \dim(L_i')=i, L_{i-1}\subset L_i'\subset L_{i+1} \} .
\end{eqnarray*}
This variety is a $\mathbb{P}^1$-bundle over the full flag variety
$\operatorname{Fl}$ in two possible ways, corresponding to
forgetting $L_i,$ respectively $L_i'.$ These two maps from
$\operatorname{Fl}_i$ onto $F$ induce two ring homomorphisms
$$ \mathrm{H}_n=\mathrm{H}(\operatorname{Fl}, \C) \lra \mathrm{H}(\operatorname{Fl}_i, \C)$$
which turn $\mathrm{H}(F_i, \C)$ into an $\mathrm{H}_n$-bimodule.
The functor $\theta_i : \mc{C}^{(1^n)} \lra \mc{C}^{(1^n)}$
is given by tensoring with this $\mathrm{H}_n$-bimodule (under the
equivalence $\mc{C}^{(1^n)}\cong \mathrm{H}_n\mathrm{-mod}$).

To describe functors $\theta_w$ for an arbitrary $w\in S_n,$ we
recall that $\operatorname{Fl}= G/B$ where $G=SL(n,\C)$ and $B$ the
Borel subgroup of $G.$ The orbits of the natural left action of $G$
on $\operatorname{Fl}\times \operatorname{Fl}$ are in natural
bijection with elements of the symmetric group. Denote by $O_w$ the
orbit associated with $w$ and by $\mathrm{IC}(\overline{O}_w)$ the
simple perverse sheaf on the closure of this orbit. The cohomology
of $\mathrm{IC}(\overline{O}_w)$ is an $\mathrm{H}_n$-bimodule, and
the functor
$$\theta_w: \mathrm{H}_n\mathrm{-mod} \lra \mathrm{H}_n\mathrm{-mod}$$
takes a module $M$ to the tensor product
$$\mathrm{H}(\mathrm{IC}(\overline{O}_w), \C)\otimes_{\mathrm{H}_n} M.$$
Notice that all cohomology rings above have a canonical grading (by cohomological
degree). The graded version of $\mc{C}^{(1^n)}$ is the category
of finite-dimensional graded $\mathrm{H}_n$-modules and the graded version of
$\theta_w$ tensors a graded module with the graded $\mathrm{H}_n$-bimodule
$\mathrm{H}(\mathrm{IC}(\overline{O}_w), \C).$

It is surprising how sophisticated the categorification of the
trivial representation is, especially when compared with the
categorification of the sign representation. Both the trivial and
the sign representation are one-dimensional, but their
categorifications have amazingly different complexities. All of the
complexity is lost when we pass to the Grothendieck group, which has
rank one.

\vspace{0.2in}

\begin{figure} [htb] \drawing{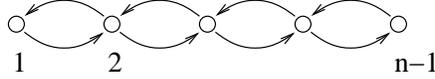} \caption{Quiver diagram
of $A_{n-1}$}\label{path2}
\end{figure}

{\bf c. Categorification of the Burau representation.} Consider the
partition $\lambda^{\ast}=(2, 1^{n-2})$ and the dual partition
$\lambda=(n-1,1).$ The category $\mc{C}^{(n-1,1)}$ admits an
explicit description, as follows. For $n>3$ let $A_{n-1}$ be the
quotient of the path algebra of the graph from Figure~\ref{path2}
by the relations
 \begin{eqnarray*}
 (i|i+1|i+2) & = & 0 , \\
 (i|i-1|i-2) & = & 0, \\
 (i|i-1|i) & = & (i|i+1|i)
\end{eqnarray*}
\vspace{0.1in}

Also, let $A_1$ be the exterior algebra on one generator, and $A_2$
be the quotient of the path algebra of the graph from Figure~\ref{path2}
(for $n=2$) by the relations $(1|2|1|2)=0=(2|1|2|1)$. The $\C$-algebra
$A_{n-1}$ is finite-dimensional.

\begin{prop} The category $\mc{C}^{(n-1,1)}$ is equivalent
to the category of finite-dimensional left $A_{n-1}$-modules.
\end{prop}

This is a well-known result, see e.g. \cite{St-quiver} for $n=2$ and
\cite{Catharina-springer} for the general case.

Denote by $P_i$ the indecomposable left projective $A_{n-1}$-module
$A_{n-1}(i).$ This module is spanned by all paths that end in vertex $i.$
Likewise, let $_iP$ stand for the indecomposable right projective $A_{n-1}$-module
$(i) A_{n-1}.$
Under the equivalence between $\mc{C}^{(n-1,1)}$ and
the category $A_{n-1}\mathrm{-mod}$ of finite-dimensional $A_{n-1}$-modules,
 the functor $\theta_i$ becomes the functor of tensoring with the bimodule
$$ P_i \otimes \hspace{0.02in} {}_iP .$$
The functors $\theta_w$ are zero for most $w\in S_n.$
They are nonzero only when the corresponding composition
of $\theta_i$'s is nonzero (which rarely happens, note that
already $\theta_i \theta_j=0$ for $|i-j|>1$).

The algebras $A_{n-1}$, as well as the modules $P_i$,  $_iP$ are
naturally graded by the length of paths. The categories of finite
dimensional graded modules over these algebras provide a
categorification of the reduced Burau representation of each of the
corresponding braid groups. For more information about the algebras
$A_{n-1}$ and their uses we refer the reader to the papers
\cite{KS}, \cite{Rou-Zimm}, \cite{Seidel-Thomas},
\cite{Catharina-corners}, \cite{HK1}.

\vspace{0.2in}


{\bf d. Categorification of the 2-column irreducible representation (partition $(2^n)$).}

\begin{figure} [htb] \drawing{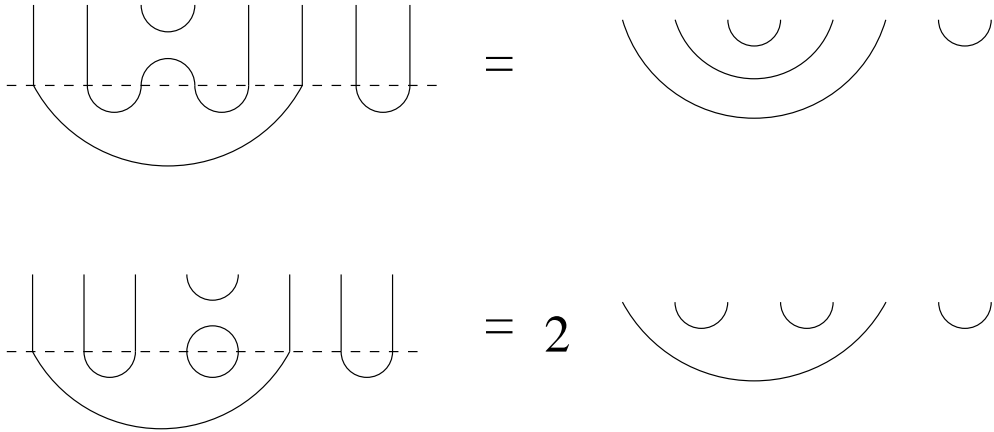} \caption{The product
$(1+s_i)b,$ for a basis element $b,$ is either another basis element
(top diagram) or the same basis element times $2.$}\label{concat} \end{figure}

Let $\lambda^{\ast}=(2^n)$ and $\lambda=(n,n).$ The irreducible
representation $S_{\Q}(\lambda^{\ast})$ has the following explicit
description. The basis of the representation consists of crossingless
matchings of $2n$ points positioned on the $x$-axis by $n$ arcs lying
in the lower half-plane, as depicted below.
\vspace{0.1in}
\begin{center}\drawing{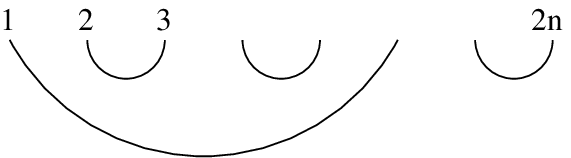} \end{center}
\vspace{0.1in}
The element $1+s_i$ acts on a basis element by concatenating
it with the diagram
\vspace{0.1in}
\begin{center}\drawing{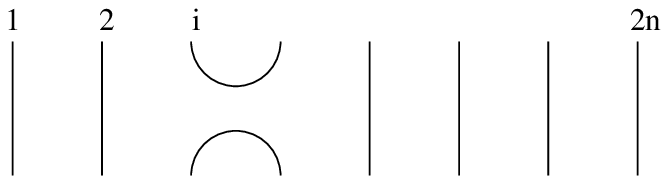} \end{center}
\vspace{0.1in}
If the concatenation contains a circle, we remove it and multiply
the result by $2,$ see figure~\ref{concat}.
\vspace{0.1in}
\vspace{0.1in}

A categorification of  this representation and of its quantum
deformation was described in \cite{Kh-FVIT}, in the context of
extending a categorification of the Jones polynomial to tangles. The
basis elements $b$ corresponding to crossingless matchings become
indecomposable projective modules $P_b$ over a certain
finite-dimensional algebra $H^n.$ The space of homs
$\Hom_{H^n}(P_a,P_b)$ between projective modules is given by gluing
crossingless matchings $a$ and $b$ along their endpoints and
applying a 2-dimensional TQFT to the resulting 1-manifold. The TQFT
is determined by a commutative Frobenius algebra, which is just the
cohomology of the 2-sphere.  Spaces of these homs together with
compositions
$$\Hom_{H^n}(P_a,P_b) \times\Hom_{H^n}(P_b,P_c)\lra \Hom_{H^n}(P_a,P_c) $$
determine $H^n$ uniquely.  The above geometric action of $1+s_i$ lifts
to the action on the category of $H^n\mathrm{-mod}$ of finite-dimensional
$H^n$-modules given by tensoring with a certain $H^n$-module.
This results in a very explicit categorification of the 2-column irreducible
representation of $S_n$ (and of the corresponding representation of the Hecke
algebra) via the category of $H^n$-modules.

It was shown in \cite{Catharina-springer} that $H^n\mathrm{-mod}$ is
equivalent to the category $\mc{C}^{(n,n)}$ generated by
projective-injective modules in the parabolic block
$\mc{O}_0^{(n,n)}$. (An extension of this equivalence to functors
will be treated in \cite{Catharina-prep}). Subquotient algebras of
$H^n$ considered  in \cite{Catharina-springer}, \cite{Chen-Kh},
\cite{Chen} can be used to categorify other 2-column representations
of the Hecke algebra and the symmetric group. These subquotient
algebras provide also a graphical description of the whole category
$\mc{O}_0^{\mu}$ for any composition $\mu_1+\mu_2=n$ of $n$
(\cite{Catharina-springer}, \cite{Catharina-prep}).


\section{Miscellaneous}

{\bf Braid group actions.}
Graded versions of projective functors $\theta_w$ categorify
the action of the Hecke algebra $H_{n,q}$ on its
various representations. There is a homomorphism
from the braid group on $n$ strands to the group
of invertible elements in $H_{n,q}.$ This homomorphism, too,
admits a categorification. The categorification should be
at least an action of the braid group on a category, and this
action is indeed well-known. To define it we need to pass
to one of the triangulated extensions of the highest weight category:
there does not seem to exist any interesting braid group actions on
abelian categories, due to the positivity imposed by the abelian
structure (see discussion in Section 6 of \cite{Kh-FVIT}).

The translation through the wall functors $\theta_i$, $i\leq 1\leq
n-1$, for the regular block $\mc{O}_0$ are compositions of two
projective endofunctors ({\it on} and {\it off the wall}) of $\cO$,
which are biadjoint to each other. This results in natural
transformations $\theta_i \lra \mathrm{Id}$ and $\mathrm{Id}\lra
\theta_i.$ Let $D(\mc{O}_0)$ be the bounded derived category of
$\mc{O}_0.$ The complexes $R_i$ and $R_i'$ of functors $0\lra
\theta_i \lra \mathrm{Id}\lra 0$ and $0 \lra \mathrm{Id}\lra
\theta_i \lra 0$ can be viewed as endofunctors of $D(\mc{O}_0)$ (we
normalise the above functors so that $\mathrm{Id}$ sits in
cohomological degree $0$).

\begin{prop} The functors $R_i$ define a braid group action on
$D(\mc{O}_0).$ The functor $R_i'$ is the inverse of $R_i.$
\end{prop}

We distinguish between \emph{weak} and \emph{genuine} group
actions; the terminology can be found in \cite{KS}, \cite{Rouquier-braid},
\cite{KT}.
That $R_i'$ and $R_i$ are inverses of each other follows from
a more general result of J.~Rickard.
That $R_i$ define a weak braid group action follows from \cite{MS-shuffling},
\cite[Proposition 10.1]{MS-root}.

The Koszul duals of the functor $R_i$ and its inverse are described
in \cite{MOS} in terms of the so-called twisting and completion
functors on $\mc{O}_0$. A geometric description of these functors
can be found in \cite{BBM-tilting} and \cite{Rouquier-braid}.

\vspace{0.2in}

The functors $\theta_i$ restrict to exact endofunctors of the
parabolic categories $\mc{O}^{\mu}_0$ and of the categories
$\mc{C}^{\mu}.$ Hence, the functors $R_i$ and $R_i'$ first induce
endofunctors on $D^b(\mc{O}^{\mu}_0),$ and define braid group
actions there and then restrict to endofunctors of the subcategory
given by complexes of projective-injective modules in
$\mc{O}^{\mu}_0$.

The braid group acts by functors respecting the triangulated
structure of the involved categories, resulting in a
categorification of parabolic braid group modules as well as those
irreducible representations of the braid group that factor through
the Hecke algebra. The two commuting actions of projective functors
on the category of Harish-Chandra bimodules as described in example
{\bf 6} of Section~\ref{sec-framework} give rise to two commuting
actions on the braid group on the derived category of the category
of Harish-Chandra bimodules. For more examples of braid group
actions on triangulated categories and a possible framework for
these actions see \cite{KT}.

\vspace{0.2in}


{\bf Invariants of tangle cobordisms.} In several cases, braid group actions on
triangulated categories can be extended to representations of the
2-category of tangle cobordisms. The objects of this 2-category (when
2-tangles are not decorated) are non-negative integers, morphisms
from $n$ to $m$ are tangles with $n$ bottom and $m$ top boundary
components, and 2-morphisms are isotopy classes of tangle cobordisms.
A representation of the 2-category of tangle cobordisms
associated a triangulated category $\mc{K}_n$ to the object $n,$
an exact functor $\mc{K}_n\lra \mc{K}_m$ to a tangle, and
a natural transformation of functors to a tangle cobordism.
Such representations can be derived from examples {\bf 7} and {\bf 8}
of Section~\ref{sec-framework}
(see \cite{Catharina-corners}, \cite{Sussan}) and from example {\bf d}
of Section~\ref{four-examples} (see \cite{Kh-FVIT}). Example {\bf 6}
 is related to at least braid
cobordisms (if not tangle cobordisms) via the construction of
\cite{Kh-Soergel}. We expect that a categorification of tensor
products of representations of quantum $\mf{sl}_2,$ mentioned at the
end of example {\bf 7} extends (after passing to derived categories,
suitable functors, and natural transformations) to a representation
of the 2-category of tangle cobordisms coloured by irreducible
representations of quantum $\mf{sl}_2.$ Such an extension would give
a categorification of the coloured Jones polynomial.

The Cautis-Kamnitzer invariant of tangle cobordisms \cite{CKa}
is based on a similar framework, but their version of the category
$\mc{K}_n$ is the derived category of coherent sheaves on a certain iterated
$\mathbb{P}^1$-bundle. The Grothendieck group of their category
is isomorphic to $V^{\otimes n},$ where $V$ is the fundamental
representation of quantum $\mf{sl}_2,$  just like in the example {\bf 7},
but these two categorifications of $V^{\otimes n}$ are noticeably
different. For instance, in the example {\bf 7} the category decomposes
into the direct sum matching the weight decomposition of the tensor product,
while the category in \cite{CKa} is indecomposable. When the
parameter is even, the two categorifications of $V^{\otimes 2n}$ appear
to have a common ``core'' subcategory, a categorification of the
invariants in $V^{\otimes 2n}$ (the latter isomorphic to $\mc{S}((2^n))$)
 briefly reviewed in the example {\bf d} above.

In the matrix factorization invariant of tangle cobordisms
\cite{KR}, the abelian category remains hidden inside the
triangulated category of matrix factorisations.

\vspace{0.2in}


{\bf Other categorifications, abelian and triangulated.} Our list of
examples of abelian categorifications is very far from complete.
Many great results in the geometric representation theory can be
interpreted as categorifications via abelian or triangulated
categories. This includes the early foundational work of
Beilinson-Bernstein and Brylinsky-Kashiwara on localisation
\cite{BB}, \cite{BK}, \cite{Mi}, the work of Kazhdan and Lusztig on
geometric realisation of representations of affine Hecke algebras
\cite{KL3}, \cite{CG}, Lusztig's geometric construction of the Borel
subalgebras of quantum groups \cite{L}, Nakajima's realisation of
irreducible Kac-Moody algebra representations as middle cohomology
groups of quiver varieties \cite{Na1}, and various constructions
related to Hilbert schemes of surfaces \cite{Grojnowski},
\cite{Na2}, quantum groups at roots of unity \cite{ABBGM}, geometric
Langlands correspondence \cite{Fr}, etc.

\vspace{0.2in}


{\bf Determinant of the Cartan matrix.} With $\lambda$ and $\mu$ as
in example {\bf 4}, let $\{P_a\}_{a\in I}$ be a collection of
indecomposable projectives in $\mc{C}^{\mu},$ one for each
isomorphism class. The Cartan matrix of $\mc{C}^{\mu}$ is an
$I\times I$ matrix $C$ with the $(a,b)$-entry being the dimension of
$\mathrm{Hom}(P_a,P_b),$ the space of homomorphisms between
projective modules $P_a$ and $P_b.$ Since $\mathrm{End}(P,P)$ is a
symmetric algebra by \cite{MS-projinj}, where $P=\oplusop{a\in
I}P_a,$ the Cartan matrix is symmetric, $c_{a,b}= c_{b,a}.$ These
algebras are not commutative, but the centre has a nice geometrical
description as the cohomology of some Springer fibre
(\cite{Kh-springer}, and more general \cite{Brundan},
\cite{Catharina-springer}).

What is the determinant of this Cartan matrix? Since $\mc{C}^{\mu}$
depends (up to equivalence) on the partition $\lambda$ only
(\cite{MS-projinj}), so does the determinant. The answer to the
question is obvious in each of the first three cases considered in
the previous section: the determinant is equal to $1$ for
$\lambda=(n),$ to $n!$ for $\lambda=(1^n),$ and to $n$ for
$\lambda=(n-1,1).$ The fourth case, when $\lambda=(2^n),$ requires
more work, and follows from the results of \cite{FGG} and
\cite{Kh-thesis}. The determinant equals
\begin{equation} \label{determinant}
  \prod_{i=1}^n  (i+1)^{r_{n,i}}, \hspace{0.2in}
   r_{n,i} = \binom{2n}{n-i}-2\binom{2n}{n-i-1}+\binom{2n}{n-i-2},
 \end{equation}
with the convention $\binom{j}{s}=0$ if $s<0.$ The answer for an
arbitrary $\lambda$ is more complicated. However, we want to point
out that this determinant of the Cartan matrix is the determinant of
the Shapovalov form (\cite{Shapovalov}) on a certain weight space of
some irreducible $\mathfrak{sl}_n$-module, as can be obtained from
instance from \cite{BK2}. It can be computed using the so-called
Jantzen-Schaper formula \cite[Satz 2]{Jakontra}.

\vspace{0.1in}

The absolute value of the determinant has an interesting categorical
interpretation. $\mc{C}^{\mu}$ is equivalent
to the category of finite-dimensional modules over some symmetric
$\C$-algebra $A^{\mu}.$  Given any symmetric $\C$-algebra
$A$ (an algebra with a nondegenerate symmetric trace $A\lra \C$),
the stable category $A\mathrm{-\underline{mod}}$ is triangulated.
Objects of $A\mathrm{-\underline{mod}}$ are finite-dimensional
$A$-modules and the set of morphisms from $M$ to $N$ is the quotient
vector space of all module maps modulo those that factor
through a projective module.  If $\det(C)\not=0$ then the Grothendieck
group of the stable category is finite abelian of cardinality
equal to the absolute value of the determinant.

\vspace{0.1in}

The graded version of this problem makes sense as well. Modules $P_a$
are naturally graded, and to a pair $(a,b)$ we can assign the
Laurent polynomial in $q$ which is the graded dimension of the graded
vector space
$\mathrm{Hom}(P_a, P_b).$ Arrange these polynomials into an
$I\times I$ matrix (the graded Cartan matrix of $\mc{C}_{gr}^{\mu}$).

\vspace{0.1in}

\emph{Problem:} Find the determinant of the graded Cartan matrix of $\mc{C}^{\mu}.$

\vspace{0.1in}

The determinant depends only on $\lambda.$
Again, the answer is known in the above four cases. In the last case,
the determinant of the graded Cartan matrix is given by formula (\ref{determinant}),
with the quantum integer $[i+1]=1+q^2+\dots + q^{2i}$ in place of
$(i+1)$ in the product (the proof follows by combining results of
\cite{FGG}  and \cite{Kh-thesis}).

The determinant is algorithmically computable, since the entries
of the graded Cartan matrix can be computed from the Kazhdan-Lusztig
polynomials of the symmetric group. We are almost tempted to conjecture
that, for any $\lambda,$ the determinant (up to a power of $q$)
is a product of quantum integers $[j]=q^{j-1}+q^{j-3}+\dots + q^{1-j},$
for small $j,$ with some multiplicities.

In \cite{JamesMathas}, a $q$-analogue of the Jantzen-Schaper formula
is obtained. Generalising \cite{Brundan} by working out a graded or
$q$-version, should imply that the determinant is equal to the
determinant of the $q$-analogue of the Shapovalov form on a suitable
weight space of an irreducible $\mf{sl}_m$-module.
 \vspace{0.2in}

{\bf Acknowledgement} {M.K. was partially supported by the NSF grant
DMS-0407784. V. M. was supported by STINT, the Royal Swedish Academy
of Sciences, and the Swedish Research Council, C. S. was partially
supported by EPSRC grant 32199. We thank all these institutions. We
also thank Kenneth Brown for useful discussions and remarks and Jon
Brundan for remarks on an earlier version of this paper.}

\vspace{0.1in}



\vspace{0.5cm}

\noindent Mikhail Khovanov, Department of Mathematics, Columbia
University, New York, NY 10027, US\\
\noindent e-mail: {\tt khovanov\symbol{64}math.columbia.edu},

\vspace{0.1in}

\noindent Volodymyr Mazorchuk, Department of Mathematics, Uppsala
University, Box 480, 751 06, Uppsala, SWEDEN, \\
\noindent e-mail: {\tt mazor\symbol{64}math.uu.se}

\vspace{0.1in}
 \noindent Catharina Stroppel, Department of
Mathematics, University of Glasgow, University Gardens, Glasgow G12
8QW, UK\\
\noindent e-mail: {\tt c.stroppel\symbol{64}maths.gla.ac.uk}

\end{document}